%% file: D_x,y_.tex
\magnification 1000
\def\arXiv{oui}
\input itenn.tex

\dimstand
\def\gs{{\got s}}
\def\gS{{\got S}}

\voffset-3mm
\input option_keys
\optionparag=1
\def\paradouze{oui}
\dateheure
\hautspages{GŽrald Tenenbaum}{Sur le biais d'une loi de probabilitŽ relative aux entiers friables}
\titrecentre{Sur le biais d'une loi de probabilitŽ relative aux entiers friables}

\bigskip\medskip
\centerline {GŽrald Tenenbaum} 
\vskip 7mm
{\eightpoint\leftskip1cm\rightskip1cm
\noi{\bf Abstract.} The standard probability law on the set $S(x,y)$ of $y$-friable integers not exceeding $x$ assigns to each friable integer $n$ a probability proportional to $1/n^\alpha$, where $\alpha=\alpha(x,y)$ is the saddle-point of the inverse Laplace integral for $\Psi(x,y):=|S(x,y)|$. This law presents a structural bias inasmuch it weights integers $>x$. We propose a quantitative measure of this bias and exhibit a related Gaussian distribution. \PAR
\medskip\noi
 \medskip
\noi{\bf Mots-clefs :} entiers friables, entiers sans grand facteur premier, mŽthode du col, Žquations diffŽrentielles aux diffŽrences, modles probabilistes. \par
\medskip
\noi \bf 2020 Mathematics Subject Classification: \rm  11N25, 11N37.\par }
\par 

\paraun{Introduction}
On dit qu'un entier naturel $n$ est $y$-friable si son plus grand facteur premier $P^+(n)$, avec la convention $P^+(1)=1$, n'excde pas $ y$.   DŽsignons par $S(x,y)$ l'ensemble des entiers $y$-friables n'excŽdant pas $x$ et par $\Psi(x,y)$ son cardinal. Posons
$$\zeta(s,y):=\prod_{p\le y}\Big(1-{1\over p^s}\Big)^{-1},\qquad \phi_y(s):=-{\zeta'(s,y)\over
\zeta(s,y)}=\sum_{p\le y}{{\log p\over p^s-1}}\qquad \big(y\ge 2,\,\re s>0\big),$$
o, ici et dans la suite, la lettre $p$ dŽsigne un nombre premier. Le point-selle
$\alpha=\alpha(x,y)$ apparaissant dans la formule asymptotique de Hildebrand--Tenenbaum~\citer{HT86}
$$\Psi(x,y)={x^\alpha\zeta(\alpha,y)\over \alpha\sqrt{2\pi|\varphi_y'(\alpha)|}}\Big\{1+O\Big({1\over
u}+{\log y\over y}\Big)\Big\}\qquad (x=y^u\ge y\ge 2)\eqdef{HT}
$$
est l'unique solution rŽelle positive de l'Žquation $\varphi_y(\alpha)=\log x$. On a  \citeplus{HT86}{(2.4)}
$$\alpha\log y=\Big\{1+O\Big({\log_22y\over \log y}\Big)\Big\}\log \Big(1+{y\over \log x}\Big)\quad(x\geqslant  y\geqslant 2).\eqdef{alphauniv}$$ 
Ici et dans la suite, nous notons $\log_k$ la $k$-ime itŽrŽe de la fonction logarithme.\par
La loi de probabilitŽ dŽfinie par 
$$\PP_{x,y}\big(\{n\}\big):=1/\{n^\alpha\zeta(\alpha,y)\}\qquad \big(P^+(n)\leqslant y\big)\eqdef{Pxy}$$
appara"t naturellement dans la thŽorie des entiers friables, notamment dans la modŽlisation des fonctions additives sur $S(x,y)$: voir  \citer{BT05a}, \citer{BT05b}, \citer{MT10}, \citer{BLT18}. La nature multiplicative de la dŽpendance en $n$ de \eqref{Pxy} est particulirement utile en pratique, mais implique que  la probabilitŽ $\PP_{x,y}$ charge des entiers $y$-friables excŽdant $x$ --- un  biais gŽnŽriquement compensŽ par la dŽcroissance rapide de $\PP_{x,y}(\{n>t\})$ lorsque $t/x\to\infty$. Une description de certaines consŽquences de ce biais et une proposition de redressement partiel apparaissent dans \citer{BT16}. 
Posant
$$D(x,y,z):=\sum_{n\in S(z,y)}{1\over n^{\alpha(x,y)}}\qquad (x\geqslant y\geqslant 2,\,z\geqslant 2),$$
nous avons ainsi $$P(x,y,z):=\PP_{x,y}\big(\{n\leqslant z\}\big)=D(x,y,z)/\zeta(\alpha,y)<1,$$ de sorte que le comportement asymptotique de cette quantitŽ, notamment lorsque $z$ est, en un sens ˆ prŽciser, voisin de $x$, fournit une mesure quantitative du biais. \par 
\par  Notre premier rŽsultat consiste en une Žvaluation de  $P(x,y,z)$ lorsque $\log (z/x)$ est convenablement dominŽ. L'ŽnoncŽ nŽcessite quelques notations supplŽmentaires.  DŽsignons par  
 $$\Phi(h):={1\over \sqrt{2\pi}}\int_{-\infty}^h\e^{-t^2/2}\d t$$
 la fonction de rŽpartition de la loi normale, notons traditionnellement $$ u:=(\log x)/\log y, \quad \ov u:=\min(y,\log x)/\log y,$$ et  posons  $$\eqalign{\Theta^2&=\Theta(x,y)^2:={|\varphi_y'(\alpha)|\over (\log y)^2}\cr&=u\Big(1+{\log x\over y}\Big)\Big\{1+O\Big({1\over \log (u+1)}+{1\over \log y}\Big)\Big\}\qquad(x\geqslant y\geqslant 2),\cr} \eqdef{defTh}$$
 o l'estimation provient de \citeplus{HT86}{(2.5)}.
DŽsignons ensuite par $\xi(t)$ la solution rŽelle positive de l'Žquation $\e^\xi=1+t\xi$ lorsque $t>1$ et convenons que $\xi(1):=0$.  La fonction $\xi$ vŽrifie
classiquement\note{Voir par exemple \citeplus{Te22}{lemme III.5.11} et la remarque qui suit l'ŽnoncŽ de \citeplus{Te22}{th. III.5.13}.}$$\xi(t)=\log t+\log_2t+O\Big({\log_2t\over \log t}\Big),\quad\xi'(t)= {1\over t}+{1\over t\log t}+O\Big({1\over t(\log t)^2}\Big)\qquad (t\to\infty).$$Notons alors que, pour tout $\varepsilon>0$ fixŽ, nous avons d'aprs \eqref{defTh} et \citeplus{HT86}{(7.19)}
$$\Theta\asymp {u\over \sqrt{\ov u}}\quad(x\geqslant y\geqslant 2),\qquad\Theta={1+O\big({1/\log y}\big)\over \sqrt{\xi'(u)}}\qquad \big((\log x)^{1+\varepsilon}<y\leqslant x\big).\eqdef{evalTh}$$
\vskip-3mm
\Propt{DZ}{Nous avons
$$\leqalignno{P(x,y,z)&\asymp 1-P(x,y,z)\asymp 1\qquad (x\geqslant y\geqslant 2, \,\log z\gg\log x,\, \log (z/x)\ll \sqrt{\ov u}\log y),&\eqdef{evalD}\cr
P(x,y,z)&=\Phi(h)+O\Big({1\over \ov u^{1/3}}\Big)\qquad (x\geqslant y\geqslant 2,\,z=xy^{\Theta h}),&\eqdef{simD+}\cr}$$
En particulier, $P(x,y,x)=\dm+O\big(1/\ov u^{1/3}\big)\quad(x\geqslant y\geqslant 2)$.
}\par \medskip\goodbreak
Nous verrons que ces rŽsultats dŽcoulent assez facilement des renseignements spŽcifiques issus de la mŽthode du col pour le comportement local  de $\Psi(x,y)$. Nous n'avons pas cherchŽ ˆ affiner le terme d'erreur de \eqref{simD+}.\note{La dŽmonstration fournit d'ailleurs une estimation lŽgrement plus prŽcise et l'on pourrait, au prix de quelques complications techniques, obtenir une dŽcroissance en $\e^{-ch^2}$ du terme d'erreur.} Cependant, la relation $\lim P(x,y,x)=\dm$ nŽcessite bien $\ov u\to\infty$. En effet, les estimations
$$\eqalign{\alpha(x,x)&=1+O\Big({1\over (\log x)^2}\Big),\quad \alpha(x,2)={1\over \log x}+O\Big({1\over (\log x)^2}\Big),\cr
\zeta(\alpha,x)&=\e^{\gamma}\log x+O(1),\qquad \zeta(\alpha,2)={\log x\over \log 2}+O(1),\cr
D(x,x,x)&=\log x +O(1),\qquad D(x,2,x)=\Big(1-{1\over \e}\Big){\log x\over \log 2}+O(1),\cr}$$
o $\gamma$ dŽsigne la constante d'Euler,
impliquent, lorsque $x\to\infty$,
$$P(x,x,x)= \e^{-\gamma}+o(1),\quad P(x,2,x)=1-{1\over \e}+o(1).$$
\par
Les Žvaluations du \ref{DZ} permettent une mesure en moyenne de la pertinence de l'approximation 
$${\Psi({x/ d},y)\over\Psi(x,y)}\approx {1\over d^\alpha}=\PP_{x,y}\big(\{n\md0d\}\big)\qquad \big(P^+(d)\leqslant y\big),$$
qui s'est rŽvŽlŽe cruciale dans nombre de problmes, comme l'estimation de la constante de Tur‡n--Kubilius friable --- voir notamment  \citer{BLT18}, \citer{BT05a}, \citer{MT10}, \citer{HMT10} --- ou encore l'Žvaluation asymptotique des valeurs moyennes friables de certaines fonctions arithmŽtiques  --- \cf, entre autres, \citer{TW03}, \citer{HTW08}, \citer{TW08a},~\citer{TW08b}. \par
\goodbreak  Posons 
$$R_d(x,y):={\Psi(x,y)\over d^\alpha}-\Psi\Big({x\over d},y\Big),\quad \sum_{d\in S(x,y)}R_d(x,y)=:\Psi(x,y)D(x,y,x)\Big\{1-\Delta(x,y)\Big\}.$$
Notant $\tau(n)$ le nombre total des diviseurs d'un entier $n$, nous avons alors
$$\Psi_\tau(x,y):=\sum_{n\in S(x,y)}\tau(n)=\sum_{d\in S(x,y)}\Psi\Big({x\over d},y\Big),\qquad \Delta(x,y)={\Psi_\tau(x,y)\over \Psi(x,y)D(x,y,x)}\cdot$$
\par 
Les estimations de Drappeau \citer{SD16} pour $\Psi_\tau(x,y)$, elles-mmes obtenues par la mŽthode du col,
permettent une estimation de $\Delta(x,y)$ valable uniformŽment pour $x\geqslant y\geqslant 2$. \par 
Nous considŽrons le domaine $$x\geqslant 3,\quad\e^{(\log_2x)^{5/3+\varepsilon}}\leqslant y\leqslant x,\leqno{(H_\varepsilon)}$$ notons  $\varrho$ la fonction de Dickman (\cf\ \citeplus{Te22}{ch. III.5}) et $\varrho_2$ son carrŽ de convolution.
\par    Le comportement asymptotique de $\varrho(t)$ est intimement liŽ ˆ celui de la fonction~$\xi(t)$ dŽfinie prŽcŽdemment.
 Notant $I(s):=\int_0^{s}(\e^v-1)\d v/v$, de sorte que (\cf\ \citeplus{Te22}{th. III.5.10})
$$\wh\varrho(-s):=\int_0^\infty\varrho(v)\e^{vs}\d v=\e^{\gamma+I(s)},$$
nous avons (\cf\ \citeplus{Te22}{th. III.5.13})
 $$\eqalign{\varrho(t)&=\Big\{1+O\Big({1\over t}\Big)\Big\}\sqrt{\xi'(t)\over 2\pi}\e^{-t\xi(t)}\wh\varrho\big(-\xi(t)\big)\cr
 &=\Big\{1+O\Big({1\over t}\Big)\Big\}\sqrt{\xi'(t)\over 2\pi}\exp\Big\{\gamma-\int_1^t\xi(s)\d s\Big\}\cr}\qquad (t\to\infty),\eqdef{asrho}$$
 une formule due ˆ Alladi \citer{Al82}, et,  d'aprs \citeplus{Sm91}{th.1 \& (3.10)} ou \citeplus{HT93}{\theoreme2}\note{En observant que $\varrho_2$ est solution de l'Žquation diffŽrentielle aux diffŽrences $t\varrho_2'(t)-\varrho_2(t)+2\varrho_2(t-1)=0$ $(t>1)$ avec la condition initiale $\varrho_2(t)=t$ $(0\leqslant t\leqslant 1)$.}, 
  $$\varrho_2(t)=\Big\{1+O\Big({1\over t}\Big)\Big\}\sqrt{\xi'(t/2)\over 4\pi}\exp\Big\{2\gamma-2\int_1^{t/2}\xi(s)\d s\Big\}=\varrho(t)2^{t+O(t/\log t)}\qquad (t\to\infty).\eqdef{asrho2}$$
DŽfinissons
$$\lambda(s,t):=\int_0^s\varrho(v)e^{v\xi(t)}\d v\qquad (s\geqslant 0,\,t\geqslant 1).$$
Nous verrons au \S\refn{prggrdy} que
$$\lambda(s,t)=\Big\{\Phi\Big(\sqrt{\xi'(t)}(s-t)\Big)+O\Big({1\over t}\Big)\Big\}\wh\varrho\big(-\xi(t)\big)\qquad (s\geqslant 0,\,t\to\infty).\eqdef{lst}$$
\par 
Posons encore
$$\eqalign{g(v)&:=v\log 4-(1+2v)\log \Big({1+2v\over 1+v}\Big)\qquad (v\geqslant 0),\cr
\vartheta(x,y)&:=\cases{\dsp 2\int_{u/2}^u\{\xi(u)-\xi(t)\}\d t&si $(x,y)\in (H_\varepsilon),$\cr
\noalign{\vskip1mm}\cr
ug\Big(\dsp{y\over \log x}\Big)& dans le cas contraire,\cr}\cr}$$
et notons que 
$$\eqalign{&\vartheta(x,y)\asymp \ov u\quad (x\geqslant y\geqslant 2),\cr
&\vartheta(x,y)=(1-\log 2)u+O(u/\log u)\quad (u\to\infty,\,y/\log x\to\infty),\cr
& \vartheta(x,y)\sim (\log 4-1)\ov u+O(\ov u^2/u)\quad (y/\log x\to0).\cr}$$
\par 
Enfin, dŽfinissons
$$\nu(t):={\varrho_2(t)\over \varrho(t)\lambda(t,t)}=\Big\{1+O\Big({1\over t}\Big)\Big\}\sqrt{2\xi'(t/2)\over \xi'(t)}
\exp\Big\{-\int_{t/2}^t\xi'(s)(2s-t)\d s\Big\}\qquad (t\geqslant 1),$$
o la formule asymptotique rŽsulte de \eqref{asrho}, \eqref{asrho2} et \eqref{lst}.\par 
\vskip-1mm
\Propt{SRd}{Posons $\varepsilon_y:=1/\sqrt{\log y}$.  Nous avons uniformŽment
$$\Delta(x,y)\asymp\e^{-\vartheta(x,y)+O(\varepsilon_y\ov u)}\qquad (x\geqslant y\geqslant 2).\eqdef{evalDelta}$$
De plus, dans le domaine $(H_\varepsilon)$, le terme $O(\varepsilon_y\ov u)$ dans \eqref{evalDelta} peut tre omis et nous avons
$$\Delta(x,y)=\nu(u)\Big\{1+O\Big({\log (u+1)\over \log y}\Big)\Big\}.\eqdef{fDeltaHep}$$} 
\par 
 La technique du prŽsent travail pourrait  tre adaptŽe ˆ l'Žtude de moyennes pondŽrŽes des restes~$R_d(x,y)$. 
\bigskip\medskip
\paraun{Preuve du \ref{DZ}}
\vskip-1mm
Lorsque $\ov u$ est bornŽ, la relation \eqref{simD+} est triviale.  Lorsque $u$ est bornŽ, \eqref{evalD} dŽcoule de
l'esti\-ma\-tion\note{Voir \eqref{alphaH} {\it infra}.} $\alpha(x,y)=1+O(1/\log x)$ et du fait que $\Psi(x,y)\asymp x$, qui dŽcoule par exemple de~\eqref{Hild} {\it infra}. Lorsque $y$ est bornŽ, nous avons $\alpha\sim \pi(y)/\log x$ et, posant $\nu_p:=1+\fl{(\log z)/\log p}$, \hbox{$\mu_p:=1+\fl{(\log z)/(\pi(y)\log p)}$}
$$\prod_{p\leqslant y}{1-p^{-\mu_p\alpha}\over 1-p^{-\alpha}}\leqslant D(x,y,z)\leqslant \prod_{p\leqslant y}{1-p^{-\nu_p\alpha}\over 1-p^{-\alpha}}\qquad \big( \log (z/x)\ll1\big),$$
ce qui implique bien \eqref{evalD}.
\par Nous pouvons donc nous placer dans la circonstance o la quantitŽ $\ov u$ est
arbitrairement grande.\par
\par 
Notant $(\log z)/\log y=u+\Theta h$,  l'hypothse relative ˆ \eqref{evalD} est $\Theta h\ll\sqrt{\ov u}$. Puisque $\Theta\asymp u/\sqrt{\ov u}$ par \eqref{defTh}, on voit que \eqref{evalD} est une consŽquence de \eqref{simD+}.
\par 
Prouvons ˆ prŽsent \eqref{simD+}. \par 
Posons $f(\sigma):=\sigma\log x+\log \zeta(\sigma,y)$, $\sigma_j:=(-1)^jf^{(j)}(\alpha)$. Nous avons
alors 
$\sigma_1=0$, $\sigma_{2}=(\Theta\log y)^2$, et (voir~\citer{HT86}), pour chaque $j\ge 2$ fixŽ,
 $$\sigma_j\asymp\ov u ^{1-j}(\log x)^j.\eqdef{sigj}$$
\par 
Notons $\alpha_v:=\alpha(y^v,y)$ et $\sigma_{j,v}=(-1)^jf^{(j)}(\alpha_v)$. La dŽrivŽe de $\alpha_v$ en fonction de $v$ vŽrifie (voir par exemple \citeplus{HT86}{(6.6)})
$$-\alpha'_v={\log y\over \sigma_{2,v}}\asymp {\ov v\over v^2\log y}\cdot\eqdef{alpha'}$$ 
\par 
Compte tenu de cette estimation, nous dŽduisons  de \citeplus{Te22}{cor. III.5.23} et \citeplus{BT05a}{th. 2.4} que 
$$\leqalignno{\Psi(t,y)/t^\alpha&\ll\Psi(x,y)/x^\alpha\qquad (t\geqslant x),&\eqdef{loc0Psi}\cr
\Psi(t,y)/t^\alpha&\asymp\Psi(x,y)/x^\alpha\qquad \big(x^{1-1/\sqrt{\ov u}}\leqslant t\leqslant x^{1+1/\sqrt{\ov u}}\,\big),&\eqdef{loc1Psi}\cr
{\Psi(t,y)\over t^\alpha}&\ll{\Psi(x,y)\over x^\alpha}\e^{-b(1-v)^2\ov u}\qquad \big(1\leqslant t=x^v\leqslant x^{1-1/\sqrt{\ov u}}\,\big),&\eqdef{loc2Psi}\cr}$$
o $b$ est une constante absolue positive. \par 
Nous pouvons clairement supposer $|h|\Theta\leqslant  u/\ov u^{1/3}$.
Il dŽcoule de \eqref{loc0Psi}, \eqref{loc1Psi}, \eqref{loc2Psi} et \eqref{HT} que
$$\eqalign{D(x,y,z)&=\int_{1-}^z{\dd\Psi(t,y)\over t^\alpha}={\Psi(z,y)\over z^\alpha}+\alpha\int_1^z{\Psi(t,y)\over t^{\alpha+1}}\d t\cr
&=\alpha\int_{x^{1-1/\ov u^{1/3}}}^z{\Psi(t,y)\over t^{\alpha+1}}\d t+O\bigg({\Psi(x,y)\over x^\alpha}+{\alpha\Psi(x,y)\log x\over x^\alpha}\int_0^{1-1/\ov u^{1/3}}\e^{-b(1-v)^2\ov u}\d v\bigg)\cr
&=\alpha\int_{x^{1-1/\ov u^{1/3}}}^z{\Psi(t,y)\over t^{\alpha+1}}\d t+O\Big({\zeta(\alpha,y)\over \sqrt{\ov u}}\Big),\cr
}$$
o nous avons utilisŽ l'estimation $\alpha\gg \ov u/(u\log y)$, qui Žquivaut ˆ $\alpha\sqrt{\sigma_2}\gg\sqrt{\ov u}$.
\par \goodbreak
Par \eqref{HT}, il suit
$$D(x,y,z)={\zeta(\alpha_u,y)\log y\over \sqrt{2\pi\sigma_{2,u}}}\int_{u-u/\ov u^{1/3}}^{u+\Theta h}{\e^{r(u,v)}}\d v+O\Big({\zeta(\alpha,y)\over \sqrt{\ov u}}\Big),\eqdef{evalD1}$$
avec $$\eqalign{r(u,v):=&\log \Big({\zeta(\alpha_v,y)\over \zeta(\alpha_u,y)}\Big)+v(\alpha_v-\alpha_u)\log y-\dm\log \Big({\sigma_{2,v}\over \sigma_{2,u}}\Big)+\log \Big({\alpha_u\over \alpha_v}\Big)\cr
=&\dm(\alpha_u-\alpha_v)^2\sigma_{2,t}+(v-u)(\alpha_v-\alpha_u)\log y+\dm\log \Big({\sigma_{2,v}\over \sigma_{2,u}}\Big)+\log \Big({\alpha_u\over \alpha_v}\Big)\cr
=&\dm(\alpha_u-\alpha_v)^2\sigma_{2,t}-(\alpha_u-\alpha_v)^2\sigma_{2,\tau}+\dm\log \Big({\sigma_{2,v}\over \sigma_{2,u}}\Big)+\log \Big({\alpha_u\over \alpha_v}\Big),\cr}$$
o $t$ et $\tau$ sont  convenablement choisis sur le segment $\gS(u,v):=\{\vartheta u+(1-\vartheta)v:0\leqslant \vartheta\leqslant 1\}$. Comme nous avons Žgalement
$$\eqalign{\max_{t\in\gS(u,v)}|\sigma_{2,t}-\sigma_{2,u}|&\ll {(u-v)\ov u\sigma_{3,u}\over u^2\log y}\ll{(u-v)\sigma_{2,u}\over u},\cr
\alpha_v-\alpha_u&={(u-v)\log y\over \sigma_{2,s}}={(u-v)\log y\over \sigma_{2,u}}\Big\{1+O\Big({u-v\over u}\Big)\Big\},\cr}$$
pour une valeur adŽquate $s\in\gS(u,v)$, il vient
$$\eqalign{r(u,v)&=-{(u-v)^2(\log y)^2\over 2 \sigma_{2,u}}\Big\{1+O\Big({1\over \ov u^{1/3}}\Big)\Big\}+O\Big({1\over \ov u^{1/3}}\Big)\cr
&=-{(u-v)^2\over 2 \Theta^2}\Big\{1+O\Big({1\over \ov u^{1/3}}\Big)\Big\}+O\Big({1\over \ov u^{1/3}}\Big),\cr}$$
o nous avons ˆ nouveau utilisŽ la minoration $\alpha_u\gg\ov u/(u\log y)$. En reportant dans \eqref{evalD1}, nous obtenons
$$D(x,y,z)=\Big\{1+O\Big({1\over \ov u^{1/3}}\Big)\Big\}\zeta(\alpha_u,y)\Phi(h)+O\Big({\zeta(\alpha,y)\over \sqrt{\ov u}}\Big)=\zeta(\alpha_u,y)\Big\{\Phi(h)+O\Big({1\over \ov u^{1/3}}\Big)\Big\}.$$
\medskip
\paraun{PrŽcisions pour les grandes valeurs de $y$}
\drefun{grdy}
La formule \eqref{simD+} peut tre affinŽe lorsque le paramtre $y$ est assez grand en fonction de $x$.
 \par 
Soit $\varepsilon>0$. Dans le domaine $(H_\varepsilon)$,
nous avons  $$\alpha=1-\xi(u)/\log y+O\big(1/(\log x\log y)\big),\eqdef{alphaH}$$ 
 \par 
Compte tenu de l'Žvaluation de Hildebrand (\cf\ \citeplus{Te22}{cor. III.5.19}) 
$$\Psi(y^v,y)=y^v\varrho(v)\Big\{1+O\Big({\log (v+1)\over \log y}\Big)\Big\}\quad((y^v,y)\in (H_\varepsilon)),\eqdef{Hild}$$ 
 et de  \citeplus{Te22}{lemme III.5.16}\note{Sous la forme $\wh{\varrho}(-\xi(u))\ll\zeta(\alpha,y)/\log y$.},  il suit, dans le domaine $(H_\varepsilon)$ pour $\log (z/x)\ll(\log x)/u^{1/3}$,
$$\eqalign{D(x,y,z)&=\alpha\int_1^z{\Psi(t,y)\over t^{\alpha+1}}\d t+O\Big({\Psi(x,y)\over x^\alpha}\Big)\cr
&=\Big\{1+O\Big({\log (u+1)\over \log y}\Big)\Big\}(\log y)\int_0^w\e^{(1-\alpha)v\log y}\varrho(v)\d v+O\Big({\zeta(\alpha,y)\over \sqrt{u}\log y}\Big)\cr&= \kappa(u,w)(\log y)\wh\varrho(-\xi(u))+O\Big({\zeta(\alpha,y)\log (u+1)\over \log y}\Big)\cr&= \zeta(\alpha,y)\Big\{\kappa(u,w)+O\Big({\log (u+1)\over \log y}\Big)\Big\},\cr}\eqdef{simD}$$
 avec toujours $w=(\log z)/\log y$ et 
 $$\kappa(u,w):={1\over \wh\varrho\big(-\xi(u)\big)}\int_0^w\e^{v\xi(u)}\varrho(v)\d v.\eqdef{defCu}$$ Pour $u=w=1$, on retrouve bien $\kappa(1,1)=\e^{-\gamma}$ puisque $\xi(1)=0$. 
\par\goodbreak 
Nous pouvons Žgalement montrer que
$$\kappa\Big(u,u+{s\over \sqrt{\xi'(u)}}\Big)=\Phi(s)+O\Big({1\over u}\Big)\qquad \big(u\geqslant 1,\,u+s/\sqrt{\xi'(u)}\geqslant 0\big).\eqdef{limC}$$
Comme le terme d'erreur de \eqref{simD} est plus prŽcis que celui de \eqref{simD+} ds que $y>\e^{(\log x)^{1/4}(\log_2x)^{3/4}}$, nous obtenons 
dans ce domaine un renforcement de \eqref{simD+}:  nous verrons que, sous cette condition, la relation \eqref{limC} implique en effet
$$\kappa(u,u+\Theta h)=\Phi(h)+O\Big({|h|\e^{-h^2/2}\over \log y}+{1\over u}\Big).\eqdef{evalkap}$$\par 
\goodbreak
Pour Žtablir \eqref{limC}, on peut supposer $s\ll \sqrt{\log 2u}$. La premire Žtape consiste ˆ majorer
$$\eqalign{E(u)&:=\int_0^{u-u^{2/3}}\e^{v\xi(u)}\varrho(v)\d v=\int_{u^{2/3}}^u\varrho(u-v)\e^{(u-v)\xi(u)}\d v\ll\int_{u^{2/3}}^u{\e^{I(\xi(u))-H(u,v)}\over \sqrt{u-v+1}}\dd v,\cr}$$
avec
$$\eqalign{H(u,v)&:=I(\xi(u))-I(\xi(u-v))-(u-v)\{\xi(u)-\xi(u-v)\}\cr
&=\int_{u-v}^u(t-u+v)\xi'(t)\d t\geqslant \dm v\log\Big( {1\over 1-v/2u}\Big)\geqslant {v^2\over 4u},\cr}$$
puisque $I(\xi(t))'=t\xi'(t)\geqslant 1$. Il suit
$$E(u)\ll\e^{I(\xi(u))}\int_{u^{2/3}}^u\e^{-v^2/4u}\d v\ll\wh\varrho\big(-\xi(u)\big)\e^{-\fs15u^{1/3}}.$$
\par \goodbreak
Pour complŽter la preuve de \eqref{limC}, il reste ˆ Žvaluer
$$\eqalign{F(u,s)&:=\int_{-s/\sqrt{\xi'(u)}}^{u^{2/3}}\e^{(u-v)\xi(u)}\varrho(u-v)\d v\cr&= \e^{\gamma+I(\xi(u))}\sqrt{\xi'(u)\over 2\pi }\int_{-s/\sqrt{\xi'(u)}}^{u^{2/3}}\e^{-H(u,v)}\Big\{1+O\Big({|v|+1\over u}\Big)\Big\}\d v,\cr}$$
puisque $\xi''(t)\ll 1/t^2$. Or
$$H(u,v)=\int_0^vw\xi'(w+u-v)\d w=\dm v^2\xi'(u)+O(v^2/u^2).$$
Il s'ensuit que
$$\eqalign{F(u,s)&= \e^{\gamma+I(\xi(u))}\sqrt{\xi'(u)\over 2\pi }\int_{-s/\sqrt{\xi'(u)}}^{u^{2/3}}\e^{-\fs12v^2\xi'(u)}\Big\{1+O\Big({|v|+1\over u}\Big)\Big\}\d v\cr&=\e^{\gamma+I(\xi(u))}\Big\{1-\Phi(-s)+O\Big({1\over u}\Big)\Big\}=\wh\varrho\big(-\xi(u)\big)\Big\{\Phi(s)+O\Big({1\over u}\Big)\Big\}.\cr}$$
Cela achve la preuve de \eqref{limC}.\par 
 La relation \eqref{evalkap} rŽsulte alors de  \eqref{evalTh}.  On notera que  $\log y\gg u^{1/3}\log (2u)$ dans le domaine considŽrŽ et qu'on peut supposer $h\ll\sqrt{\log 2u}$.
\medskip\medskip
\paraun{Preuve du \ref{SRd}} 
D'aprs \citeplus{SD16}{th.1}, nous avons
$$\Psi_\tau(x,y)=\beta\sqrt{\pi\gs_2}\Psi\big(\sqrt{x},y\big)^2\Big\{1+O\Big({1\over \ov u}\Big)\Big\}\qquad (x\geqslant y\geqslant 2)$$
o l'on a posŽ $\beta:=\alpha\big(\sqrt{x},y\big)$, $\gs_2:=|\varphi'_y(\beta)|$. Compte tenu de \eqref{HT} et \eqref{evalD}, il suit
$$\Delta(x,y)\asymp {x^{\beta-\alpha}\zeta(\beta,y)^2\over \zeta(\alpha,y)^2},$$
o nous avons fait usage de l'estimation $\alpha\sqrt{\sigma_2}\asymp\beta\sqrt{\gs_2}$, qui dŽcoule de \citeplus{HT86}{th. 2}.
\par \goodbreak
Posons
$$\eqalign{&\vartheta_1(x,y):=\normalbaselineskip=18pt\cases{\xi(u)-\xi(u/2),& si $(x,y)\in(H_\varepsilon)$,\cr
\dsp\log\Big(1+{y\over y+\log x}\Big),& dans le cas contraire,}
\cr
&\vartheta_2(x,y):=\normalbaselineskip=20pt\cases{\dsp2\int_{u/2}^ut\xi'(t)\d t,&si $x\geqslant y> (\log x)(\log_22x)^3$,\cr
\dsp{2y\over \log y}\log \Big(1+{\log x\over 2y+\log x}\Big),& dans le cas contraire,}\cr
&\vartheta_0(x,y):=
\vartheta_2(x,y)-u\vartheta_1(x,y).\cr}$$
Nous avons alors, avec la notation $\varepsilon_y$ de l'ŽnoncŽ,
$$\vartheta_0(x,y)=\normalbaselineskip=18pt\cases{\vartheta(x,y) & si $(x,y)\in (H_\varepsilon)$ ou $y\leqslant (\log x)(\log_2x)^3$,\cr
\vartheta(x,y)+O\big(u \varepsilon_y\big)& dans le cas contraire.\cr}$$
\par 
La premire Žtape de la preuve de \eqref{evalDelta} consiste ˆ montrer que
$$(\beta-\alpha)\log x=u\vartheta_1(x,y)+O(1+\varepsilon_y\ov u)\qquad (x\geqslant y\geqslant 2),\eqdef{b-a}$$
o le terme $O(\varepsilon_y\ov u)$ peut tre omis si $(x,y)\in(H_\varepsilon)$.\par 
Conservons la notation $\alpha_v:=\alpha(y^v,y)$, de sorte que $\alpha=\alpha_u$, $\beta=\alpha_{u/2}$. Dans le domaine $(H_\varepsilon)$, il rŽsulte de \eqref{alphaH} que
$$\beta-\alpha={\xi(u)-\xi(u/2)\over \log y}+O\Big({1\over \log x}\Big).$$
Lorsque, disons, $\log y\leqslant (\log_2x)^2$, l'Žvaluation 
$$\sigma_{2,v}=\Big(1+{v\log y\over y}\Big)v(\log y)^2\Big\{1+O\Big({1\over \log (v+1)}\Big)\Big\}$$
Žtablie en \citeplus{HT86}{th. 2} fournit, compte tenu de \eqref{alpha'},
$$\eqalign{(\beta-\alpha)\log y&=\int_{u/2}^u{(\log y)^2\over \sigma_{2,v}}\d v=\{1+O(\varepsilon_y)\}y\int_{u/2}^u{\dd v\over v(y+v\log y)}\cr
&=\{1+O(\varepsilon_y)\}\log \Big(1+{y\over y+\log x}\Big).\cr}$$
Cela implique bien \eqref{b-a}.\par 
Il reste ˆ Žtablir l'estimation
$$Z(x,y):=\log\Big( {\zeta(\alpha,y)\over \zeta(\beta,y)}\Big)=\dm\vartheta_2(x,y)+O(1+\varepsilon_y\ov u)\qquad (x\geqslant y\geqslant 2),\eqdef{Z}$$
avec la mme convention que pour \eqref{b-a} concernant le terme d'erreur $O(\varepsilon_y\ov u)$.
\par \goodbreak
Dans le domaine $(H_\varepsilon)$, l'estimation \eqref{Z} rŽsulte du calcul effectuŽ dans \citeplus{HT86}{\page289}. Dans le domaine complŽmentaire, l'Žvaluation de $\varphi_y'(\sigma)$ apparaissant au lemme 13 de \citer{HT86} permet d'Žcrire
$$Z(x,y)=\Big\{1+O\Big({1\over \log y}\Big)\Big\}\int_\alpha^\beta{y^{1-\sigma}-1\over (1-\sigma)(1-y^{-\sigma})}\d\sigma.\eqdef{Z1}$$
La preuve donnŽe dans \citer{HT86} de l'Žvaluation \eqref{alphaH} pour $\alpha=\alpha_u$ fournit en fait (\cf\ \citeplus{HT86}{pp. 285-7}), pour une constante $C$ assez grande,
$$\alpha(x,y)=1-{\xi(u)\over \log y}+O\Big({1\over (\log y)^2}\Big)\qquad \big(C(\log x)(\log_2x)^3< y\leqslant x\big).\eqdef{alpha+}$$
Cela implique dans le mme domaine $y^\alpha\gg y\e^{-\xi(u)}\gg\log y$, de sorte que le terme $1-y^{-\sigma}$ peut tre englobŽ par le terme d'erreur. Le changement de variables $(1-\sigma)\log y=\xi(t)$ et l'Žvaluation~\eqref{alpha+} fournissent alors \eqref{Z}.\par 
Lorsque $y\leqslant C(\log x)(\log_22x)^3$, nous effectuons le changement de variables $t=y^\sigma$ dans l'intŽgrale de~\eqref{Z1}. Il vient
$$\eqalign{Z(x,y)&=\Big\{1+O\Big({1\over \log y}\Big)\Big\}\int_{y^\alpha}^{y^\beta}{y-t\over t(t-1)\log (y/t)}\d t.\cr}$$
Or l'Žvaluation \eqref{alphauniv}  implique $\beta\log y\ll\log_22y $, d'o
$$\eqalign{Z(x,y)=\Big\{1+O\Big({\log_22y\over \log y}\Big)\Big\}{y\over \log y}\int_{y^\alpha}^{y^\beta}{\dd t\over t(t-1)}=\Big\{1+O\Big({\log_22y\over \log y}\Big)\Big\}{y\over \log y}\log \Big({1-y^{-\beta}\over 1-y^{-\alpha}}\Big).\cr}$$
Comme nous avons \citeplus{HT86}{(7.8)}
$${1\over 1-y^{-\alpha_v}}=\Big\{1+O\Big({1\over \log y}\Big)\Big\}{y+v\log y\over y}\qquad (v\geqslant 1,\,y\geqslant 2),$$
cela fournit ˆ nouveau \eqref{Z}. 
\par 
Pour achever la preuve du \ref{SRd},  observons que l'Žvaluation \eqref{fDeltaHep} rŽsulte de l'avant-dernire formule \eqref{simD} avec $w=u$, de \eqref{limC} avec $s=0$, de la formule de Hildebrand \eqref{Hild}, et de la formule analogue pour $\Psi_\tau(x,y)$, telle qu'elle dŽcoule de \citeplus{TW03}{cor. 2.3}.
\medskip\medskip
\paraun{Remarque finale} 
\vskip-1mm
Il rŽsulte de \eqref{evalD}, \eqref{simD} et \eqref{limC} que, pour une constante absolue convenable $c>0$,
$$c\leqslant \kappa(u,u)\leqslant 1-c\qquad (u\geqslant 1).$$
 On peut Žtablir cela directement, sans invoquer \eqref{limC}.
En effet, d'aprs \citeplus{FT96}{lemme 6.1} et \hbox{\citeplus{Te22}{\thinspace(III.5.63)}}, nous avons, pour une constante absolue convenable $c_0$,
$$\varrho(u-v)\asymp \varrho(u)\e^{v\xi(u)}\quad(1\leqslant v\leqslant c_0\sqrt{u}).$$
Cela implique d'une part
$$\int_0^u\e^{v\xi(u)}\d v\gg\e^{u\xi(u)}\int_0^{c_0\sqrt{u}}\e^{-v\xi(u)}\varrho(u-v)\d v\asymp\varrho(u)\e^{u\xi(u)}\sqrt{u}\asymp \wh{\varrho}(-\xi(u))$$
et d'autre part
$$\eqalign{\int_u^\infty\varrho(v)\e^{v\xi(u)}\d v&=\e^{u\xi(u)}\int_0^\infty\varrho(u+v)\e^{v\xi(u)}\d v\cr
&\gg\varrho(u)\e^{u\xi(u)}\int_0^{c_0\sqrt{u}}\e^{v\xi(u)-v\xi(u+v)}\d v\gg\varrho(u)\e^{u\xi(u)}\sqrt{u}\asymp\wh{\varrho}(-\xi(u)).\cr}$$
\par 
\noi{\bf Remerciements.} L'auteur tient ici ˆ exprimer ici sa plus  plus chaleureuse gratitude envers   l'arbitre pour sa relecture attentive et ses suggestions pertinentes.
\goodbreak
\bigskip\medskip
\centerline{\twelvebf Bibliographie}
\bigskip
{\eightpoint\leftskip9mm\rightskip5mm
\bibtem{Al82} K. Alladi, The Tur\'an--Kubilius inequality for integers without large prime factor,
{\it J. reine angew. Math. \bf 335} (1982), 180--196.
\bibtem{BLT18} R. de la Bretche, Y. Lamzouri \& G. Tenenbaum, InŽgalitŽ de Tur‡n-Kubilius friable et indŽpendance asymptotique, \AA\ \bf183, \rm \numero2 (2018), 191-199.
\bibtem{BT05a} R. de la Bretche \& G. Tenenbaum, PropriŽtŽs statistiques des entiers friables, {\it Ramanujan J. \bf9} (2005), 139--202.
\bibtem{BT05b} R. de la Bretche \& G. Tenenbaum, Entiers friables : inŽgalitŽ de Tur‡n--Kubilius et 
applications (avec R. de la Bretche), {\it Invent. Math. \bf159} (2005), 531--588.
\bibtem{BT16} R. de la Bretche \& G. Tenenbaum, Sur l'inŽgalitŽ de Tur‡n--Kubilius friable, {\it J. London Math. Soc.} (2) {\bf 93}, \numero 1 (2016), 175--193.
\bibtem{SD16} S. Drappeau, Remarques sur les moyennes des fonctions de Piltz sur les entiers friables; {\it Quart. J. Math. (Oxford), \bf 67}, \numero4 (2016), 507-517.
\bibtem{FT96} ƒ. Fouvry \& G. Tenenbaum, RŽpartition statistique des entiers sans grand facteur premier
dans les progressions arithmŽtiques,
{\it Proc. London Math. Soc.}, (3) \bf72 \rm(1996) 481-514.
\bibtem{HMT10}  G. Hanrot, B. Martin \& G. Tenenbaum,  Constantes de Tur‡n-Kubilius friables : Žtude numŽrique, {\it Exp. Math.} Experiment. Math. \bf19, \numero3 \rm(2010),  345Ð361.
\bibtem{HTW08}  G. Hanrot, G. Tenenbaum \& J. Wu, Moyennes de certaines fonctions
multiplicatives sur les entiers friables, 2, \PLMS\ (3) \bf96 \rm(2008) 107--135.
\bibtem{HT86}  A. Hildebrand \& G. Tenenbaum, On integers free of large prime factors, {\it
Trans. Amer. Math. Soc. \bf 296} (1986), 265--290.
\bibtem{HT93}  A. Hildebrand \& G. Tenenbaum, On a class of difference differential equations arising in number theory, {\it J.~d'Analyse \bf 61} (1993),  145--179.
\bibtem{MT10} B. Martin \& G. Tenenbaum, Sur l'inŽgalitŽ de
Tur‡n-Kubilius friable, {\it J. reine angew. Math. \bf647} (2010), 175Ð234.
\bibtem{Sm91} H. Smida, Sur les puissances de convolution de la fonction de Dickman,
\AA\ {\bf 59}, \numero2  (1991), 124--143.
\bibtem{Te22} G. Tenenbaum, {\it Introduction ˆ la thŽorie analytique et probabiliste des
nombres}, 5me Žd., Dunod, 2022.\par 
\bibtem{TW03}G. Tenenbaum \& J. Wu, Moyennes de certaines fonctions
multiplicatives sur les entiers friables, {\it J. reine angew. Math. \bf 564} (2003), 119--166.
\bibtem{TW08a} G. Tenenbaum \& J. Wu, Moyennes de certaines fonctions
multiplicatives sur les entiers friables, 3, {\it Compositio Math. \bf144} \numero2 (2008), 339--376.
\bibtem{TW08b} G. Tenenbaum \& J. Wu, Moyennes de certaines fonctions
multiplicatives sur les entiers friables, 4, {\it Actes du colloque de MontrŽal, 2006}, Centre de Recherches MathŽmatiques, CRM
Proceedings and Lecture Notes 46 (2008), 129--141.
\par}
\vskip 5mm
{\sevenrm\baselineskip9pt
GŽrald Tenenbaum\par
Institut ƒlie Cartan\par 
UniversitŽ de Lorraine\par
 BP 70239\par
54506 VandÏuvre-ls-Nancy Cedex\par
 France
\smallskip
\seventt gerald.tenenbaum@univ-lorraine.fr\par}

\bye

%% file: itenn.tex
  
 \def\oui{oui} 
  
\ifx\arXiv\oui
\else
 \pdfpagewidth=210truemm
 \pdfpageheight=297truemm 
\fi
  
  %
  %

  %
  %

  \catcode`@=12 

 \def\defrefnote#1{\definexref{#1}{{\the\footnotenumber}}{refnotes}}

  %
  %


\ifx\couleurs\oui
\input graphicx
\fi

\input eplain.tex
 
\ifx\couleurs\oui
\beginpackages
\usepackage{color}
\endpackages 
\long\def\rge#1{{\color{red}#1}}

\definecolor{bleu-iecn}{cmyk}{.98,.13,.1,.55}

\fi

\makeatletter
\def\numberedfootnote{%
ÊÊ\global\advance\footnotenumber by 1
ÊÊ\@eplainfootnote{{\number\footnotenumber}}%
}%
\def\makecolumns#1/#2 {\par \begingroup
ÊÊ \@columndepth = #1
ÊÊ \advance\@columndepth by -1
ÊÊ \divide \@columndepth by #2
ÊÊ \advance\@columndepth by 1
ÊÊ \@linestogoincolumn = \@columndepth
ÊÊ \@linestogo = #1
ÊÊ \currentcolumn = 1
ÊÊ \def\@endcolumnactions{%
ÊÊÊÊÊÊ\ifnum \@linestogo<2
ÊÊÊÊÊÊÊÊ \the\crtok \egroup \endgroup \par 
ÊÊÊÊÊÊ\else
ÊÊÊÊÊÊÊÊ \global\advance\@linestogo by -1
ÊÊÊÊÊÊÊÊ \ifnum\@linestogoincolumn<2
ÊÊÊÊÊÊÊÊÊÊÊÊ\global\advance\currentcolumn by 1
ÊÊÊÊÊÊÊÊÊÊÊÊ\global\@linestogoincolumn = \@columndepth
ÊÊÊÊÊÊÊÊÊÊÊÊ\the\crtok
ÊÊÊÊÊÊÊÊ \else
ÊÊÊÊÊÊÊÊÊÊÊÊ&\global\advance\@linestogoincolumn by -1
ÊÊÊÊÊÊÊÊ \fi
ÊÊÊÊÊÊ\fi
ÊÊ }%
ÊÊ \makeactive\^^M
ÊÊ \letreturn \@endcolumnactions
ÊÊ \@columnwidth = \hsize
ÊÊÊÊ \advance\@columnwidth by -\parindent
ÊÊÊÊ \divide\@columnwidth by #2
ÊÊ \penalty\abovecolumnspenalty
ÊÊ \noindent 
ÊÊ \valign\bgroup
ÊÊÊÊ &\hbox to \@columnwidth{\strut \hsize = \@columnwidth ##\hfil}\cr
}%
\makeatother

\lefteqnumbers
   \def\testd{oui}
   \def\choixlat{\ifx\numadroite\testd\righteqnumbers
            \else  \lefteqnumbers\fi}
    \choixlat

\catcode`@=\letter
\def\@eplainfootnote#1{\let\@sf\empty 
  \ifhmode\edef\@sf{\spacefactor\the\spacefactor}\/\fi
  \global\advance\hlfootlabelnumber by 1
  \hlstart@impl{foot}{\hlfootlabel}%
  \hldest@impl{footback}{\hlfootbacklabel}%
  \hbox{$^{(#1)}$}%
  \hlend@impl{foot}%
  \@sf\vfootnote{#1.}%
}%
\catcode`@=\other

  \interfootnoteskip=0pt
  \let\note=\numberedfootnote
  \everyfootnote={\eightpoint\leftskip=5truemm\rightskip5truemm}
  
  \hsize150truemm\vsize 240truemm\hoffset=5truemm
  
  \def\dimstand{\hsize 150truemm\vsize 250truemm\hoffset=8truemm\voffset=5truemm}
  
  \pretolerance=500\tolerance=1000\brokenpenalty=5000
  \parindent3mm
  
  \countdef\temps=170
  \temps=\time
  \countdef\nminutes=171{\nminutes = \time}
  \countdef\nheure=172
  \def\heure{\begingroup                     
     \temps = \time \divide\temps by 60
     \nheure = \temps                        
     \nminutes = \time
     \multiply\temps by 60
     \advance\nminutes by -\temps            
     \ifnum\nminutes<10 \toks1 = {0}%
     \else\toks1 = {}%
     \fi
     \number\nheure h\the\toks1 \number\nminutes  
  \endgroup}%

  \newcount\chstart
  \chstart=\pageno
 \headline={\ifnum\pageno=\chstart {\hfill} \else {\hss \tenrm --\ \folio\ --\hss}\fi}
  \footline={\hfill}
  \normalbaselines
  \frenchspacing
    \def\dater{\vglue-10mm\rightline{(\the\day/\the\month/\the\year)}}
  \def\dateheure{\vglue-10mm\rightline{(\the\day/\the\month/\the\year,\ \heure)}}

  \newif\ifpagetitre \pagetitretrue
\newtoks\hautpagetitre \hautpagetitre={\hfill}
\newtoks\baspagetitre \baspagetitre={\hfill}
\newtoks\auteurcourant \auteurcourant={\hfill}
\newtoks\titrecourant \titrecourant={\hfill}
\newtoks\hautpagegauche
\newtoks\hautpagedroite
\newtoks\hautpagemilieu
\hautpagemilieu={\tenrm\hfil -- \folio\ -- \hfil}
\hautpagegauche={\ifx\midfolio\oui\the\hautpagemilieu\else\tenrm\folio\hfill\the\auteurcourant\hfill\fi}
\hautpagedroite={\ifx\midfolio\oui\the\hautpagemilieu\else\hfill\the\titrecourant\hfill\tenrm\folio\fi}
\newtoks\baspagegauche \baspagegauche={\hfil}
\newtoks\baspagedroite \baspagedroite={\hfil}
\headline={\ifpagetitre\the\hautpagetitre
\else\ifodd\pageno\the\hautpagedroite\else\the\hautpagegauche\fi\fi }
\footline={\ifpagetitre\the\baspagetitre
\else\ifodd\pageno\the\baspagedroite
\else\the\baspagegauche\fi\fi \global\pagetitrefalse}

\def\pageblanche{\vfill\eject\pagetitretrue
\null\vfill\eject
\pagetitretrue
}
\def\chgtpage{\ifodd\pageno \else
\pageblanche \fi \pagetitretrue\titreun=0\footnotenumber=0}

\def\chgtpageincrtitreun{\ifodd\pageno \else
\pageblanche \fi \pagetitretrue\footnotenumber=0}

\def\majnombres{\ifodd\pageno \else
\pageblanche \fi \pagetitretrue\hautpoly\titreun=0\footnotenumber=0}

\def\hautspages#1#2{\auteurcourant={\ninepcap#1}\titrecourant={\nineit#2}}

\ifnum\chstart=\pageno \pagetitretrue\fi
  

  \def\PAR{\par}
  
  \def\leftnote#1{\vadjust{\setbox1=\vtop{\hsize 20mm\parindent=0pt\eightpoint
  \baselineskip=9pt\rightskip=4mm plus 4mm\vskip-4.7mm#1}\hbox{\kern-2cm\smash{\box1}}}}

  
  \def\raggedcenter{\leftskip=20pt plus 10em  
       \rightskip=\leftskip 
        \parfillskip=0pt 
         \spaceskip=.3333em \xspaceskip=.5em 
          \pretolerance=9999 \tolerance=9999
           \hyphenpenalty=9999 \exhyphenpenalty=9999 }
           
  \def\titrecentre#1{{\parindent0mm\raggedcenter
       \spaceskip=.6em plus .2em minus .2em\xspaceskip=.6em plus .2em minus .2em
        \tit#1\par}}
        


  \def\oui{oui}
  
\def\fontetitreun{\ifx\paradouze\oui\douzepts\gpdouze\twelvebf\textfont1=\twelveib\else
\quatorzepts\gpquatorze\fourteenbf\fi}

\def\fontetitreunl{\douzepts\textfont1=\twelveib\scriptfont1=\tenib\fourteenti}
 
 \def\fontetitredeux{\textfont1=\eleveni\ifx\paradouze\oui\onzepts\scriptfont1=\ninei\elevenit\else
                        \douzepts\twelveit\fi}
 
   \def\fontetitredeuxb{\ifx\paradouze\oui\onzepts\eleventi\gponze\textfont1=\elevenib\scriptfont1=\nineib
                         \else\douzepts\twelveti\scriptfont1=\twelveib\scriptfont1=\tenib\gpdouze\fi}
                         
\def\fontetitredeuxl{\onzepts\textfont1=\elevenbf\scriptfont1=\ninebf\twelvebf}
  
\def\fontetitretrois{\textfont0=\elevenrm\scriptfont0=\eightrm\textfont1=\eleveni
                      \scriptfont1=\eighti\scriptscriptfont1=\sixi\elevenit}
                      
\def\fontetitrequatre{\textfont0=\elevenrm\scriptfont0=\eightrm\textfont1=\eleveni
                      \scriptfont1=\eighti\scriptscriptfont1=\sixi\elevenrm}
  
  \newcount\titreun\titreun=0
  \newcount\titredeux\titredeux=0
  \newcount\titretrois\titretrois=0
  \newcount\titrequatre\titrequatre=0
  \newcount\enonce\enonce=0
  
  \def\incr#1{\global\advance#1 by 1 {\the #1}}
  \def\avance#1{\global\advance#1 by 1}
  \def\init#1{\global#1=0}
  
  \long\def\Indentation#1#2{\setbox10=\hbox{\fontetitreun#1}
  	                    \ifdim\wd10 < 4mm
                         \setbox10=\hbox to 4mm{\box10\hfill}
                       \else\ifdim\wd10 < 6mm
                         \setbox10=\hbox to 6mm{\box10\hfill}
  	                    \else\ifdim\wd10 < 8mm
                         \setbox10=\hbox to 8mm{\box10\hfill}
                       \else\ifdim\wd10 < 12mm
                         \setbox10=\hbox to 12mm{\box10\hfill}
                       \fi\fi\fi\fi
                       \dimen10=\hsize
                       \advance \dimen10 by -\wd10
                       \noindent \box10 %
                       \ignorespaces
                       \hbox{\vtop{\hsize=\dimen10\raggedright\noindent\fontetitreun#2}}}

  \long\def\paraun#1{\removelastskip\par\medskip\goodbreak\vskip0pt plus.01\vsize\penalty-100
                \vskip0pt plus-.01\vsize
  	              \init{\titredeux}\ifnum\optionparag=1{\init\eqnumber\init\enonce}\else{}\fi
                  \goodbreak{\fontetitreun
  	                \Indentation{\incr{\titreun}.\ }{\fontetitreun #1\par}}\nobreak\medskip}

 %
 %
 \long\def\paraunc#1{\removelastskip\par\bigskip\goodbreak\vskip0pt plus.01\vsize\penalty-100
                \vskip0pt plus-.01\vsize
  	              \init{\titredeux}
                 \ifnum\optionparag=1{\init{\eqnumber}\init\enonce}\else{}\fi
                  \goodbreak
  	                {\parindent0mm\raggedcenter\fontetitreun\incr{\titreun}.\ 
                     \fontetitreun #1\par}\nobreak\medskip}
                     
\newtoks\titreunl
\titreunl={\ifnum\titreun=1{I}\fi%
\ifnum\titreun=2{II}\fi%
\ifnum\titreun=3{III}\fi%
\ifnum\titreun=4{IV}\fi%
\ifnum\titreun=5{V}\fi%
\ifnum\titreun=6{VI}\fi%
\ifnum\titreun=7{VII}\fi%
\ifnum\titreun=8{VIII}\fi%
\ifnum\titreun=9{IX}\fi%
\ifnum\titreun=10{X}\fi%
\ifnum\titreun=11{XI}\fi%
\ifnum\titreun=12{XII}\fi%
\ifnum\titreun=13{XIII}\fi%
}
\long\def\paraunl#1{\removelastskip\par\bigskip\bigskip\goodbreak\vskip0pt plus.01\vsize\penalty-100
                \vskip0pt plus-.01\vsize
  	              \init{\titredeux}\ifnum\optionparag=1{\init\eqnumber\init\enonce}\else{}\fi
                  \goodbreak{\fontetitreunl
  	                \Indentation{\global\advance\titreun by 1{\the\titreunl}.\ }{\fontetitreunl #1\par}}\nobreak\smallskip}

  
  \long\def\paradeux#1{\init{\titretrois}\vskip0pt plus.01\vsize\penalty-10
                \vskip0pt plus-.01\vsize\ifx \elie\oui\medskip\ifnum\titredeux>0\medskip\fi\fi
                 \Indentation{\fontetitredeux\the\titreun${\cdot}$\incr{\titredeux}.}
                              {\fontetitredeux\textfont1=\eleveni#1}\nobreak\par }
  
  \long\def\paradeuxb#1{\init{\titretrois}\vskip0pt plus.001\vsize\penalty-10
                \vskip0pt plus-.01\vsize{\ifx \elie\oui\medskip\ifnum\titredeux>0\medskip\fi\fi
                  \Indentation
  {\fontetitredeuxb\the\titreun${\cdot}$\incr{\titredeux}.}{ \fontetitredeuxb#1}}\nobreak
\smallskip}

\newtoks\titredeuxl
\titredeuxl={\ifnum\titredeux=1{A}\fi%
\ifnum\titredeux=2{B}\fi%
\ifnum\titredeux=3{C}\fi%
\ifnum\titredeux=4{D}\fi%
\ifnum\titredeux=5{E}\fi%
\ifnum\titredeux=6{F}\fi%
\ifnum\titredeux=7{G}\fi%
\ifnum\titredeux=8{H}\fi%
\ifnum\titredeux=9{I}\fi%
\ifnum\titredeux=10{J}\fi%
\ifnum\titredeux=11{K}\fi%
\ifnum\titredeux=12{L}\fi%
\ifnum\titredeux=13{M}\fi%
}
 \long\def\paradeuxl#1{\init{\titretrois}\vskip0pt plus.001\vsize\penalty-10
                \vskip0pt plus-.01
                \vsize \bigskip%
                  \Indentation
     {\fontetitredeuxl\global\advance\titredeux by 1
  \quad \the\titreunl${\cdot}$\the\titredeuxl.}{ \fontetitredeuxl#1}
  \removelastskip\nobreak\smallskip}
  

  \long\def\paratrois#1{\init{\titrequatre}\ifdim\lastskip<\smallskipamount
                \removelastskip\smallskip\fi
                 \vskip0pt plus.01\vsize\penalty-10
                  \vskip0pt
plus-.01\vsize{\ifx \elie\oui\ifnum\titretrois>0\medskip\fi\fi
\Indentation{\fontetitretrois\the\titreun${\cdot}$\the\titredeux${\cdot}$\incr{\titretrois}.\ }
  {\hskip0mm\baselineskip=14pt\fontetitretrois#1}\nobreak\smallskip}}
  
  
  \long\def\paratroisl#1{\init{\titrequatre}\ifdim\lastskip<\smallskipamount
                \removelastskip\fi
                 \vskip0pt plus.01\vsize\penalty-10
                  \vskip0pt
plus-.01\vsize\ifx \elie\oui\bigskip
\fi
\Indentation{\fontetitretrois\quad \quad \the\titreunl{${\cdot}$}\the\titredeuxl${\cdot}$\incr{\titretrois}.\ }
  {\hskip0mm\fontetitretrois#1}\nobreak\smallskip}


  \long\def\paraquatre#1{\ifdim\lastskip<\smallskipamount
                \removelastskip\smallskip\fi
                 \vskip0pt plus.01\vsize\penalty-10
                  \vskip0pt
                  plus-.01\vsize\par
 
\Indentation{\fontetitrequatre \the\titreun{${\cdot}$}\the\titredeux${\cdot}$\the\titretrois${\cdot}$\incr{\titrequatre}.\ }
{\hskip0mm\fontetitrequatre#1}\nobreak\smallskip}


\newtoks\titrequatrel
\titrequatrel={\ifnum\titrequatre=1{a}\fi%
\ifnum\titrequatre=2{b}\fi%
\ifnum\titrequatre=3{c}\fi%
\ifnum\titrequatre=4{d}\fi%
\ifnum\titrequatre=5{e}\fi%
\ifnum\titrequatre=6{f}\fi%
\ifnum\titrequatre=7{g}\fi%
\ifnum\titrequatre=8{h}\fi%
\ifnum\titrequatre=9{i}\fi%
\ifnum\titrequatre=10{j}\fi%
\ifnum\titrequatre=11{k}\fi%
\ifnum\titrequatre=12{l}\fi%
\ifnum\titrequatre=13{m}\fi%
}
\long\def\paraquatrel#1{\ifdim\lastskip<\smallskipamount
                \removelastskip\smallskip\fi
                 \vskip0pt plus.01\vsize\penalty-10
                  \vskip0pt
                  plus-.01\vsize{\bigskip
\Indentation{\global\advance\titrequatre by 1
\fontetitrequatre\quad \quad \quad \the\titreunl${\cdot}$\the\titredeuxl${\cdot}$\the\titretrois${\cdot}$\the\titrequatrel.\ }
{\hskip0mm\fontetitrequatre#1}\nobreak\smallskip}}

\ifx\optionkeys\oui
\def\drefun#1{\definexref{¤#1}{{\the\titreun}}{}} 
\def\drefdeux#1{\definexref{¤#1}{{\the\titreun}.{\the\titredeux}}{}}
\def\dreftrois#1{\definexref{¤#1}{{\the\titreun}.{\the\titredeux}.{\the\titretrois}}{}}
\else
\def\drefun#1{\definexref{prg#1}{{\the\titreun}}{}} 
\def\drefdeux#1{\definexref{prg#1}{{\the\titreun}.{\the\titredeux}}{}}
\def\dreftrois#1{\definexref{prg#1}{{\the\titreun}.{\the\titredeux}.{\the\titretrois}}{}}
\fi

%


  \long\def\propdeux#1#2#3#4{%
       \avance{\enonce}
       \leavevmode\edef\temp{#2}%
         \ifx\temp\empty 
          \else
           \definexref{#2}{#1~{\the\titreun.\the\enonce}}{enonces}
            \definexref{s#2}{{\the\titreun.\the\enonce}}{enonces}
             \fi
\smallskip
      \noindent{\bf#1\ {\bf\the\titreun.\the\enonce{#3}.}\enspace}{\sl#4\par}%
      \ifdim\lastskip<\medskipamount \removelastskip\penalty55\par \fi
   }

  \long\def\propun#1#2#3#4{%
      \avance{\enonce}
       \leavevmode\edef\temp{#2}%
        \ifx\temp\empty 
          \else
           \definexref{#2}{#1~{\the\enonce}}{enonces}
            \definexref{{s#2}}{{\the\enonce}}{enonces}
             \fi
   \par 
     \noindent{\bf#1\ {\bf\the\enonce{#3}.}\enspace}{\sl#4\par}%
     \ifdim\lastskip<\medskipamount \removelastskip\penalty55\medskip\fi
  }
  
  \long\def\prop#1#2#3#4{\ifnum\optionparag=1
                          \propdeux{#1}{#2}{\textfont1=\elevenib#3}{#4} \else\propun{#1}{#2}{\textfont1=\elevenib#3}{#4}\fi}

  \long\def\propt#1#2#3{\ifx\tpf\oui \prop{Th\'eo\-r\`eme}{#1}{#2}{#3}\par
                       \else\prop{Theorem}{#1}{#2}{#3}\par\fi}
  \long\def\Propt#1#2{\propt{#1}{}{#2}}
  \long\def\propl#1#2#3{\ifx\tpf\oui\prop{Lem\-me}{#1}{#2}{#3}\par
                         \else\prop{Lemma}{#1}{#2}{#3}\par\fi}
  
  \long\def\propc#1#2#3{\ifx\tpf\oui\prop{Corol\-laire}{#1}{#2}{#3}\par
                         \else\prop{Corollary}{#1}{#2}{#3}\par\fi}

  \long\def\propd#1#2#3{\ifx\tpf\oui\prop{D\'efi\-nition}{#1}{#2}{#3}\par
                       \else\prop{Definition}{#1}{#2}{#3}\par\fi} 
  
  \long\def\proptd#1#2#3{\ifx\tpf\oui\prop{Th\'eor\`eme et d\'efi\-nition}{#1}{#2}{#3}\par
                       \else\prop{Theorem and definition}{#1}{#2}{#3}\par\fi}


  
  \newcount\optionparag\optionparag=1
  
  \long\def\section#1#2{\ifnum\optionparag=1 \paraun{#2} 
                        \else\goodbreak{\fontetitreun
  	                \Indentation{#1.\ }{#2}}\nobreak\smallskip\fi}

  \def\eqconstruct#1{\ifnum\optionparag=1{\the\titreun\hbox{$\cdot$}#1}\else{#1}\fi}

  
  
  \def\numref{oui}  
  
  \newcount\mesref\mesref=0 
  \def\defbib#1{\ifx\numref\oui\global\advance\mesref by 1 \definexref{#1}{{\the
                 \mesref}}{}\else\definexref{#1}{#1}{}\fi}
  \def\bibtem#1{\defbib{#1}\item{\citer{#1}}}
  \def\citer#1{[\ref{#1}]}
  \def\citeplus#1#2{[\ref{#1}; #2]}

  
  \font\seventeenmsa=msam10 at 17pt    
  \font\fourteenmsa=msam10 at 14pt
  \font\twelvemsa=msam10 at 12pt
  \font\tenmsa=msam10                 
  \font\ninemsa=msam10 at 9pt 
  \font\eightmsa=msam10 at 8pt 
  \font\sevenmsa=msam7 
  \font\sixmsa=msam10 at 6pt
  \font\fivemsa=msam5
  \newfam\msafam\textfont\msafam=\tenmsa\scriptfont\msafam=\sevenmsa\scriptscriptfont\msafam=\fivemsa
  
  \font\seventeenbb=msbm10 at 17pt     
  \font\fourteenbb=msbm10 at 14pt
  \font\twelvebb=msbm10 at 12pt
  \font\tenbb=msbm10                   
  \font\ninebb=msbm10 at 9pt 
  \font\eightbb=msbm10 at 8pt 
  \font\sevenbb=msbm7 
  \font\sixbb=msbm10 at 6pt
  \font\fivebb=msbm5 
  \newfam\bbfam\textfont\bbfam=\tenbb\scriptfont\bbfam=\sevenbb\scriptscriptfont\bbfam=\fivebb
  \def\bb{\fam\bbfam\tenbb}%

  \font\seventeenscaln=eusm10 at 17pt   
  \font\twelvescaln=eusm10 at 12pt
  \font\tenscaln=eusm10                
  \font\ninescaln=eusm10 scaled 900
  \font\eightscaln=eusm10 scaled 800
  \font\sevenscaln=eusm10 scaled 700
  \font\sixscaln=eusm10 scaled 600
   
  \newfam\scalnfam\textfont\scalnfam=\tenscaln\scriptfont\scalnfam=\sevenscaln\scriptscriptfont\scalnfam=\sixscaln
  \def\scaln{\fam\scalnfam\tenscaln}%
  \def\scal{\scaln}
  
  \font\tenscalb=eusb10                

  \font\sevenscalb=eusb10 scaled 700

  \newfam\scalbfam\textfont\scalbfam=\tenscalb\scriptfont\scalbfam=\sevenscalb
  %
  
  %
  %
  \font\fourteenrm=cmr12 scaled 1200
  \font\elevenrm=cmr10 at 11pt
  \font\twelverm=cmr12
  \font\ninerm=cmr9
  \font\eightrm=cmr8      
  \font\sevenrm=cmr7
  \font\sixrm=cmr6

  \font\seventeenpcap=cmcsc10 at 17pt
  \font\tenpcap=cmcsc10                        
  \font\ninepcap=cmcsc9
  \font\eightpcap=cmcsc8
  \font\sevenpcap=cmcsc10 scaled 700
  
  \newfam\pcapfam\textfont\pcapfam=\tenpcap\scriptfont\pcapfam=\sevenpcap
  \def\pcap{\fam\pcapfam\tenpcap}
  
  \font\seventeenrm=cmbx12 scaled 1400

  \font\fourteenbf=cmbx10 scaled 1400
  
  \font\twelvebf=cmbx12
  \font\elevenbf=cmbx10 at 11pt
  \font\ninebf=cmbx9  
  \font\eightbf=cmbx8
  \font\sixbf=cmbx6
  
  \font\tengot=eufm10                           
   
  \font\eightgot=eufm10 at 8truept 
  \font\sevengot=eufm7 
  \font\sixgot=eufm10 at 6 truept 
   
  \newfam\gotfam
  \textfont\gotfam=\tengot\scriptfont\gotfam=\sevengot\scriptscriptfont\gotfam=\sixgot
  \def\got{\fam\gotfam\tengot}%

  
  \def\tit{%
  \textfont0=\seventeenrm\scriptfont0=\tenrm\def\rm{\fam0\seventeenrm}%
  \textfont1=\seventeenib\scriptfont1=\twelveib%
  \textfont2=\seventeensy\scriptfont2=\twelvesy\scriptscriptfont2=\ninesy
  \textfont3=\seventeenex
  \textfont\itfam=\seventeenti
  \def\it{\fam\itfam\seventeenti}%
  \textfont\bbfam=\seventeenbb \scriptfont\bbfam=\twelvebb
  \def\bb{\fam\bbfam\seventeenbb}%
  \textfont\msafam=\seventeenmsa\scriptfont\msafam=\twelvemsa
  \textfont\scalnfam=\seventeenscaln
  \def\pcap{\fam\pcapfam\seventeenpcap}
  \normalbaselineskip=25pt\normalbaselines\rm}

  \font\seventeenti=cmbxti10 scaled 1680
  
  \font\fourteenti=cmbxti10 at 14pt
  
  \font\twelveti=cmbxti10 scaled 1200
  \font\eleventi=cmbxti10 at 11pt

  %
  %
  \font\twelveit=cmti12	
  \font\elevenit=cmti10 scaled 1100
  \font\nineit=cmti9
  \font\eightit=cmti8
  \font\sevenit=cmti7

  %
  %
 
 \font\seventeenib=cmmib10 scaled 1680
  \font\fourteenib=cmmib10 scaled 1400
  \font\twelveib=cmmib10 scaled 1200
  \font\elevenib=cmmib10 scaled 1100
  \font\tenib=cmmib10
\font\eightib=cmmib10 scaled 800
  \font\nineib=cmmib10 scaled 900
\font\sevenib=cmmib10 scaled 700
\font\sixib=cmmib10 scaled 600
\font\fiveib=cmmib10 scaled 500

\ifx\ITAN\oui
\else
\innernewfam\cmmibfam
\textfont\cmmibfam=\tenib
\scriptfont\cmmibfam=\sevenib
\scriptscriptfont\cmmibfam=\fiveib
\def\ib{\fam\cmmibfam\tenib}
\fi

  %
  %
  \font\twelvei=cmmi10 scaled 1200
  \font\eleveni=cmmi10 scaled 1100
  \font\ninei=cmmi9
  \font\eighti=cmmi8 
  \font\seveni=cmmi7 			                
  \font\sixi=cmmi6
  
  \font\ninesl=cmsl9                    
  \font\eightsl=cmsl8 
  \font\sevensl=cmsl10 at 7pt

  \font\ninett=cmtt9                    
  \font\eighttt=cmtt8
  \font\seventt=cmtt10 scaled 700

  \font\seventeensy=cmsy10 scaled 1680    
  \font\fourteensy=cmsy10 scaled 1400
  \font\twelvesy=cmsy10 scaled 1176
  
  \font\ninesy=cmsy9                      
  \font\eightsy=cmsy8
  \font\sixsy=cmsy6
  \font\seventeenex=cmex10 at 17pt
  \font\fourteenex=cmex10 at 14pt
  \font\twelveex=cmex10 at 12pt
  \font\nineex=cmex10 at 9pt
  \font\eightex=cmex10 at 8pt
  \font\sevenex=cmex10 at 7pt
  \font\sixex=cmex10 at 6pt
  \font\fiveex=cmex10 at 5pt
  
   
  \font\fourteengp=cmmi10 at 14pt
  
  \font\twelvegp=cmmib10 at 12pt
  \font\elevengp=cmmib10 at 11pt
  \font\tengp=cmmib10                          
  \font\ninegp=cmmib10 at 9pt 
  \font\eightgp=cmmib8 
   
  \font\sixgp=cmmib6


  \def\gponze{\textfont0=\elevenbf\scriptfont0=\eightbf\scriptscriptfont0=\sixbf
           \textfont1=\elevengp\scriptfont1=\eightgp\scriptscriptfont1=\sixgp}
  \def\gpdouze{\textfont0=\twelvebf\scriptfont0=\tenbf\scriptscriptfont0=\ninebf
           \textfont1=\twelvegp\scriptfont1=\tengp\scriptscriptfont1=\ninegp}        
  
 \def\gpquatorze{\textfont0=\fourteenbf\scriptfont0=\twelvebf\scriptscriptfont0=\elevenbf
           \textfont1=\fourteengp\scriptfont1=\twelvegp\scriptscriptfont1=\elevengp}

  
  \expandafter\chardef\csname pre amssym.def at\endcsname=\the\catcode`\@
  \catcode`\@=11
  \def\undefine#1{\let#1\undefined}
  \def\newsymbol#1#2#3#4#5{\let\next@\relax
   \ifnum#2=\@ne\let\next@\msafam@\else
   \ifnum#2=\tw@\let\next@\bbfam@\fi\fi
   \mathchardef#1="#3\next@#4#5}
  \def\mathhexbox@#1#2#3{\relax
   \ifmmode\mathpalette{}{\m@th\mathchar"#1#2#3}%
   \else\leavevmode\hbox{$\m@th\mathchar"#1#2#3$}\fi}
  \def\hexnumber@#1{\ifcase#1 0\or 1\or 2\or 3\or 4\or 5\or 6\or 7\or 8\or
   9\or A\or B\or C\or D\or E\or F\fi}
  
  \def\setboxz@h{\setbox\z@\hbox}
  \def\wdz@{\wd\z@}
  \def\boxz@{\box\z@}
  
  \edef\msafam@{\hexnumber@\msafam}
  \mathchardef\dabar@"0\msafam@39
  
  \edef\bbfam@{\hexnumber@\bbfam}
  \def\widehat#1{\setboxz@h{$\m@th#1$}%
   \ifdim\wdz@>\tw@ em\mathaccent"0\bbfam@5B{#1}%
   \else\mathaccent"0362{#1}\fi}
  \def\widetilde#1{\setboxz@h{$\m@th#1$}%
   \ifdim\wdz@>\tw@ em\mathaccent"0\bbfam@5D{#1}%
   \else\mathaccent"0365{#1}\fi}
  \newsymbol\leqq 1335          
  \newsymbol\leqslant 1336
  \newsymbol\lessgtr 1337       
  \newsymbol\backprime 1038     
  \newsymbol\risingdotseq 133A  
  \newsymbol\fallingdotseq 133B 
  \newsymbol\succcurlyeq 133C   
  \newsymbol\geqq 133D          
  \newsymbol\geqslant 133E
  \newsymbol\nmid 232D
  \newsymbol\nexists 2040
  \newsymbol\smallsetminus 2272
  \newsymbol\varnothing 203F 
  \catcode`\@=\active

  \catcode`\@=11
  \newcount\typofr\typofr=1
  
  \catcode`\;=\active
  \def;{\ifnum\typofr=1\relax\ifhmode\ifdim\lastskip>\z@\unskip\fi
     \kern.2em\fi\string;\else\string;\fi}
  
  \catcode`\:=\active
  \def:{\ifnum\typofr=1\relax\ifhmode\ifdim\lastskip>\z@\unskip\fi
  \penalty\@M\ \fi\string:\else\string:\fi}
  
  \catcode`\!=\active
  \def!{\ifnum\typofr=1\relax\ifhmode\ifdim\lastskip>\z@\unskip\fi
     \kern.2em\fi\string!\else\string!\fi}
  
  \catcode`\?=\active
  \def?{\ifnum\typofr=1\relax\ifhmode\ifdim\lastskip>\z@\unskip\fi
     \kern.2em\fi\string?\else\string?\fi}

  \def\francais{\typofr=1\def\tpf{oui}}
  
  \def\oui{oui}
  \francais
  
  \catcode`\@=12
  

\ifx\textures\oui
\def\raye #1|{\leavevmode\setbox1=\hbox{#1}%
\raise .5pt\hbox to \wd1{\xleaders\hbox{\rge{$ \sct / $}%
\kern 1pt}\hfill\kern -1pt }\kern -\wd1 \unhbox1\relax }

\def\barre #1|{\leavevmode\setbox1=\hbox{#1}%
\rlap{\Red\vrule height 2.4pt depth -1.2pt width \wd1}\Black \unhbox1\relax}
\else
\def\raye #1|{\leavevmode\setbox1=\hbox{#1}%
\raise .5pt\hbox to \wd1{\xleaders\hbox{\rge{$ \sct / $}%
\kern 1pt}\hfill\kern -1pt }\kern -\wd1 \unhbox1\relax }

\def\barre #1|{\leavevmode\setbox1=\hbox{#1}%
\rlap{\color{red}\vrule height 2.4pt depth -1.2pt width \wd1}\color{black} \unhbox1\relax}

\fi
  

  
  \def\og{\leavevmode\raise.24ex\hbox{$\scriptscriptstyle\langle\!\langle\>$}}    
  \def\fg{\leavevmode\raise.24ex\hbox{$\scriptscriptstyle\>\rangle\!\rangle$}}    

  \def\d{\,{\rm d}}
  \def\dd{{\rm d}}

  \def\PP{{\bb P}}

  \def\O{{\scal O}}
  \def\P{{\scaln P}}

  \def\frac#1#2{{#1\over #2}}

  \def\numero{n$^{\rm o}\thinspace$}

  \def\page#1{\rm p.\thinspace#1}
  \def\theoreme#1{\rm th.\thinspace#1}
  

  \def\numero{n$^{\rm o}\thinspace$}

  \def\¤{\S\thinspace}

  \def\¥{$\bullet$ }
  
  
  \def\e{{\rm e}}
  
  \def\md#1#2{\equiv#1\,({\rm mod\,}#2)}

  \def\epsilon{\varepsilon}

  \def\phi{\varphi}
  \def\theta{\vartheta}
  \def\rho{\varrho}
  \def\dm{{\textstyle{1\over 2}}}
  
  \def\dsp{\displaystyle}
  \def\sct{\scriptstyle}
  \def\pf{\noi{\it Proof. }}
  \def\nid{\ifnum\typofr=1\par\noindent{\it D\'emonstration. }\else\pf\fi}
  \def\noi{\noindent}
  \def\rem{\ifnum\typofr=1\noi{\it Remarque.}\ \else\noi{\it Remark.}\ \fi}
  \def\rems{\ifnum\typofr=1\noi{\it Remarques.}\ \else\noi{\it Remarks.}\ \fi}
  \def\re{{\Re e\,}}
  
  \def\ov{\overline}

  \def\1{{\bf 1}}
  \def\|{\Vert}

  \def\le{\leqslant}\def\leq{\leqslant}
  \def\ge{\geqslant}\def\geq{\geqslant}
  \def\wh{\widehat}
  \def\cf{{cf.}}

  \newsymbol\subsetneqq 2324
  \newsymbol\subsetneq 2328

  \def\fl#1{\left\lfloor #1 \right\rfloor}

  \def\log{\mathop{\rm log}\nolimits}
  
  \def\fs#1#2{{\scriptstyle{#1\over #2}}}
  



  \def\pmb#1{\setbox0=\hbox{#1}%
  \kern-.025em\copy0\kern-\wd0\kern.05em\copy0\kern-\wd0\kern-.025em\raise .0433em\box0 }

  
  \skewchar\eighti='177 \skewchar\sixi='177
  \skewchar\eightsy='60 \skewchar\sixsy='60
  
  \def\eightpoint{%
  \textfont0=\eightrm\scriptfont0=\sixrm\scriptscriptfont0=\fiverm
  \def\rm{\fam0\eightrm}%
  \textfont1=\eighti\scriptfont1=\sixi
  \scriptscriptfont1=\fivei\def\oldstyle{\fam1\seveni}%
  \textfont2=\eightsy\scriptfont2=\sixsy\scriptscriptfont2=\fivesy
  \textfont3=\eightex\scriptfont3=\sixex
  \textfont\itfam=\eightit
  \def\it{\fam\itfam\eightit}%
  \textfont\slfam=\eightsl
  \def\sl{\fam\slfam\eightsl}%
  \textfont\bbfam=\eightbb \scriptfont\bbfam=\sixbb\scriptscriptfont\bbfam=\fivebb
  \def\bb{\fam\bbfam\eightbb}%
  \textfont\msafam=\eightmsa\scriptfont\msafam=\sixmsa
  \textfont\scalnfam=\eightscaln
  \def\scaln{\fam\scalnfam\eightscaln}
  \textfont\ttfam=\eighttt
  \def\tt{\fam\ttfam\eighttt}%
\textfont\gotfam=\eightgot
  \textfont\bffam=\eightbf\scriptfont\bffam=\sixbf\scriptscriptfont\bffam=\fivebf
  \def\bf{\fam\bffam\eightbf}%
  \ifx\ITAN\oui\else\textfont\cmmibfam=\eightib
       \scriptfont\cmmibfam=\sixib
        \scriptscriptfont\cmmibfam=\fiveib
         \def\ib{\fam\cmmibfam\eightib}
   \fi
  \textfont\pcapfam=\eightpcap
  \def\pcap{\fam\pcapfam\eightpcap}
  \abovedisplayskip=2pt plus2pt minus 2pt
  \belowdisplayskip=2pt plus1pt minus 2pt
  \abovedisplayshortskip= 1pt plus 2pt minus 1pt
  \belowdisplayshortskip= 1pt plus 2pt minus 1pt
  \smallskipamount=2pt plus 1pt minus 2pt
  \medskipamount=3pt plus 2pt minus 2pt
  \bigskipamount=7pt plus 3pt minus 3pt
  \setbox\strutbox=\hbox{\vrule height 5pt depth 2pt width 0pt}%
  \normalbaselineskip=9pt\normalbaselines\rm}

  \def\({\left(}
  \def\){\right)}
  
  \def\footnoterule{\kern -2pt\hrule width 7truecm\kern 2.4pt}
  
  \def\xnotedef#1{\definexref{#1}{\noexpand\number\footnotenumber}{Note}}%

  
  
  \def\ninepoint{%
  \textfont0=\ninerm\scriptfont0=\sixrm\scriptscriptfont0=\fiverm
  \def\rm{\fam0\ninerm}%
  \textfont1=\ninei\scriptfont1=\sixi
  \scriptscriptfont1=\fivei\def\oldstyle{\fam1\ninei}%
  \textfont2=\ninesy\scriptfont2=\sixsy\scriptscriptfont2=\fivesy
  \textfont3=\nineex\scriptfont3=\sixex
  \textfont\itfam=\nineit
  \def\it{\fam\itfam\nineit}%
  \textfont\slfam=\ninesl
  \def\sl{\fam\slfam\ninesl}%
  \textfont\bbfam=\ninebb\scriptfont\bbfam=\sixbb\scriptscriptfont\bbfam=\fivebb
  \def\bb{\fam\bbfam\ninebb}%
  \textfont\msafam=\ninemsa\scriptfont\msafam=\sixmsa\scriptscriptfont\msafam=\fivemsa
  \textfont\scalnfam=\ninescaln
  \def\scaln{\fam\scalnfam\ninescaln}
  \textfont\ttfam=\ninett
  \def\tt{\fam\ttfam\ninett}%
  \textfont\bffam=\ninebf\scriptfont\bffam=\sixbf\scriptscriptfont\bffam=\fivebf
  \def\bf{\fam\bffam\ninebf}%
  \abovedisplayskip=3pt plus2pt minus 2pt
  \belowdisplayskip=3pt plus1pt minus 2pt
  \abovedisplayshortskip= 2pt plus 2pt minus 1pt
  \belowdisplayshortskip= 2pt plus 2pt minus 1pt
  \smallskipamount=2pt plus 1pt minus 2pt
  \medskipamount=3pt plus 2pt minus 2pt
  \bigskipamount=7pt plus 3pt minus 3pt
  \setbox\strutbox=\hbox{\vrule height 5pt depth 2pt width 0pt}%
  \normalbaselineskip=11pt plus.3pt minus.2pt\normalbaselines\rm}

  \def\sevenpoint{%
  \textfont0=\sevenrm\scriptfont0=\sixrm\scriptscriptfont0=\fiverm
  \def\rm{\fam0\sevenrm}%
  \textfont1=\seveni\scriptfont1=\sixi
  \scriptscriptfont1=\fivei\def\oldstyle{\fam1\seveni}%
  \textfont2=\sevensy\scriptfont2=\sixsy\scriptscriptfont2=\fivesy
  \textfont3=\sevenex\scriptfont3=\fiveex
  \textfont\itfam=\sevenit
  \def\it{\fam\itfam\sevenit}%
  \textfont\slfam=\sevensl
  \def\sl{\fam\slfam\sevensl}%
  \textfont\bbfam=\sevenbb \scriptfont\bbfam=\sixbb\scriptscriptfont\bbfam=\fivebb
  \def\bb{\fam\bbfam\sevenbb}%
  \textfont\msafam=\sevenmsa\scriptfont\msafam=\sixmsa
  \textfont\scalnfam=\sevenscaln
  \def\scaln{\fam\scalnfam\sevenscaln}
  \textfont\bffam=\sevenbf\scriptfont\bffam=\sixbf\scriptscriptfont\bffam=\fivebf
  \def\bf{\fam\bffam\sevenbf}%
  \textfont\ttfam=\seventt
  \abovedisplayskip=2pt plus2pt minus 2pt
  \belowdisplayskip=2pt plus1pt minus 2pt
  \abovedisplayshortskip= 1pt plus 2pt minus 1pt
  \belowdisplayshortskip= 1pt plus 2pt minus 1pt
  \smallskipamount=2pt plus 1pt minus 2pt
  \medskipamount=3pt plus 2pt minus 2pt
  \bigskipamount=7pt plus 3pt minus 3pt
  \setbox\strutbox=\hbox{\vrule height 5pt depth 2pt width 0pt}%
  \normalbaselineskip=9pt\normalbaselines\rm}

 \def\onzepts{%
 \textfont0=\elevenrm\scriptfont0=\ninerm
 \textfont1=\eleveni\scriptfont1=\ninei
}

\def\douzepts{%
  \textfont0=\twelverm\scriptfont0=\tenrm\def\rm{\fam0\twelverm}%
  \textfont1=\twelvei\scriptfont1=\teni%
  \textfont2=\twelvesy\scriptfont2=\tensy\scriptscriptfont2=\eightsy
  \textfont3=\twelveex
  \textfont\itfam=\twelveti
  \def\it{\fam\itfam\twelveti}%
  \textfont\bffam=\twelvebf\scriptfont\bffam=\tenbf\scriptscriptfont\bffam=\eightbf
  \def\bf{\fam\bffam\twelvebf}%
  \textfont\bbfam=\twelvebb \scriptfont\bbfam=\tenbb
  \def\bb{\fam\bbfam\twelvebb}%
  \textfont\msafam=\twelvemsa\scriptfont\msafam=\tenmsa
  \textfont\scalnfam=\twelvescaln
  \normalbaselineskip=15pt\normalbaselines\rm}

\def\quatorzepts{%
  \textfont0=\fourteenrm\scriptfont0=\twelverm\def\rm{\fam0\fourteenrm}%
  \textfont1=\fourteenib\scriptfont1=\twelveib%
  \textfont2=\fourteensy\scriptfont2=\twelvesy\scriptscriptfont2=\tensy
  \textfont3=\fourteenex
  \textfont\itfam=\fourteenti
  \def\it{\fam\itfam\fourteenti}%
  \textfont\bffam=\fourteenbf\scriptfont\bffam=\twelvebf\scriptscriptfont\bffam=\tenbf
  \def\bf{\fam\bffam\fourteenbf}%
  \textfont\bbfam=\fourteenbb \scriptfont\bbfam=\twelvebb
  \def\bb{\fam\bbfam\fourteenbb}%
  \textfont\msafam=\fourteenmsa\scriptfont\msafam=\twelvemsa
  \textfont\scalnfam=\twelvescaln
  \normalbaselineskip=18pt\normalbaselines\rm}


\def\AA{{\it Acta Arith.}}

\def\PLMS{{\it Proc. London Math. Soc.}}

\def\picture #1 by #2 (#3){\leavevmode\vbox to #2{
     \hrule width #1 height 0pt depth 0pt
      \vfill
       \special{picture #3}}}

\def\illustration #1 by #2 (#3) scaled #4{\dimen1=#2
  \divide\dimen1 by 1000
  \multiply\dimen1 by #4
  \vtop to \dimen1{\dimen1=#1
  \divide\dimen1 by 1000
  \multiply\dimen1 by #4
  \hsize=\dimen1\vss
  \noindent\special{illustration #3 scaled #4}}}

\ifx\couleurs\oui

\fi

%% file: option_keys.tex
\ifx\optionkeymacros\undefined\else \fi

\catcode`\Œ=\active\defŒ{{\aa}}       
\catcode`\º=\active\defº{\int}        
\catcode`\=\active\def{\c c}        
\catcode`\¶=\active\def¶{\partial}    
\catcode`\Ä=\active\defÄ{\oint}       
\catcode`\Æ=\active\defÆ{\triangle}   
\catcode`\Â=\active\defÂ{\neg}        
\catcode`\µ=\active\defµ{\mu}         
\catcode`\¿=\active\def¿{{\o}}        
\catcode`\¹=\active\def¹{\pi}         
\catcode`\Ï=\active\defÏ{{\oe}}       
\catcode`\§=\active\def§{{\ss}}       
\catcode`\ =\active\def {\dagger}     
\catcode`\Ã=\active\defÃ{\sqrt}       
\catcode`\·=\active\def·{\Sigma}      
\catcode`\Å=\active\defÅ{\approx}     
\catcode`\½=\active\def½{\Omega}      
\catcode`\£=\active\def£{{\it\$}}     
\catcode`\°=\active\def°{\infty}      
\catcode`\¤=\active\def¤{{\S}}        
\catcode`\¦=\active\def¦{{\P}}        
\catcode`\¥=\active\def¥{\bullet}     
\catcode`\»=\active\def»{\leavevmode\raise.585ex\hbox{\b a}}      
\catcode`\¼=\active\def¼{\leavevmode\raise.6ex\hbox{\b o}}        
\catcode`\­=\active\def­{\not=}       
\catcode`\²=\active\def²{\leq}        
\catcode`\³=\active\def³{\geq}        
\catcode`\Ö=\active\defÖ{\div}        
\catcode`\É=\active\defÉ{{\dots}}     
\catcode`\¾=\active\def¾{{\ae}}       
\catcode`\Ç=\active\defÇ{\og}         
\catcode`\Ò=\active\defÒ{``}          
\catcode`\Á=\active\defÁ{!`}          
\catcode`\¢=\active\def¢{\rlap/c}     
\catcode`\Ô=\active\defÔ{`}           
\catcode`\Õ=\active\defÕ{'}           


\catcode`\=\active\def{{\AA}}       
\catcode`\'=\active\def'{\c C}        
\catcode`\¯=\active\def¯{{\O}}        
\catcode`\¸=\active\def¸{\Pi}         
\catcode`\Î=\active\defÎ{{\OE}}       
\catcode`\®=\active\def®{{\AE}}       
\catcode`\×=\active\def×{\diamond}    
\catcode`\¡=\active\def¡{\accent'27}  
\catcode`\Ó=\active\defÓ{''}          
\catcode`\±=\active\def±{\pm}         
\catcode`\È=\active\defÈ{\fg}         
\catcode`\À=\active\defÀ{?`}          
\catcode`\Ð=\active\defÐ{--}          
\catcode`\Ñ=\active\defÑ{---}         


\catcode`\Š=\active\defŠ{\"a}        
\catcode`\'=\active\def'{\"e}        
\catcode`\•=\active\def•{\"{\i}}     
\catcode`\š=\active\defš{\"o}        
\catcode`\Ÿ=\active\defŸ{\"u}        
\catcode`\Ø=\active\defØ{\"y}        
\catcode`\å=\active\defå{\^A}        
\catcode`\€=\active\def€{\"A}        
\catcode`\…=\active\def…{\"O}        
\catcode`\†=\active\def†{\"U}        
\catcode`\‡=\active\def‡{\'a}        
\catcode`\Ž=\active\defŽ{\'e}        
\catcode`\'=\active\def'{\'{\i}}     
\catcode`\—=\active\def—{\'o}        
\catcode`\œ=\active\defœ{\'u}        
\catcode`\ƒ=\active\defƒ{\'E}        
\catcode`\æ=\active\defæ{\^E}        
\catcode`\é=\active\defé{\`E}        %
\catcode`\ˆ=\active\defˆ{\`a}        
\catcode`\=\active\def{\`e}        
\catcode`\"=\active\def"{\`{\i}}     
\catcode`\˜=\active\def˜{\`o}        
\catcode`\=\active\def{\`u}        
\catcode`\Ë=\active\defË{\`A}        
\catcode`\‹=\active\def‹{\~a}        
\catcode`\–=\active\def–{\~n}        
\catcode`\›=\active\def›{\~o}        
\catcode`\Ì=\active\defÌ{\~A}        
\catcode`\"=\active\def"{\~N}        
\catcode`\Í=\active\defÍ{\~O}        
\catcode`\‰=\active\def‰{\^a}        
\catcode`\=\active\def{\^e}        
\catcode`\"=\active\def"{\^{\i}}     
\catcode`\™=\active\def™{\^o}        
\catcode`\ž=\active\defž{\^u}        

\let\optionkeymacros\null